\documentclass[a4paper, 10pt, titlepage, roman]{article}
\usepackage{MRHStyle}
\usepackage{JournalOfExpMath}

\usepackage[square]{natbib}
\setcitestyle{citesep={;}}

\renewcommand{\phi}{\ensuremath{\varphi}}

\newcommand{\rhomin}{\ensuremath{\rho_{min}}}
\newcommand{\rhomax}{\ensuremath{\rho_{max}}}
\newcommand{\logf}[1]{\ensuremath{\log\inps{#1}}}

\begin{document}
\title{\Large \textbf{A constructive proof of the convergence of Kalantari’s bound on polynomial zeros}}
\author{\Large Matt Hohertz \\ Department of Mathematics, Rutgers University \\ mrh163@math.rutgers.edu}

\maketitle

\begin{abstract}
In his 2006 paper, Jin proves that Kalantari's bounds on polynomial zeros, indexed by $m\geq 2$ and called $L_m$ and $U_m$ respectively, become sharp as $m\rightarrow\infty$.  That is, given a degree $n$ polynomial $p(z)$ not vanishing at the origin and an error tolerance $\epsilon > 0$, Jin proves that there exists an $m$ such that
	\begin{equation*}
		\frac{L_m}{\rhomin} \geq 1 - \epsilon,
	\end{equation*}
where $\rhomin:=\min_{\rho : p(\rho) = 0} \aval{\rho}$.  In this paper we derive a formula that yields such an $m$, thereby constructively proving Jin's theorem.  In fact, we prove the stronger theorem that this convergence is uniform in a sense, its rate depending only on $n$ and a few other parameters.   We also give experimental results that suggest an optimal $m$ of (asymptotically) $O\inps{\frac{1}{\epsilon^d}}$ for some $d\ll 2$.  A proof of these results would show that Jin's method runs in $O\inps{\frac{n}{\epsilon^d}}$ time, making it efficient for isolating high-degree polynomial zeros.
\end{abstract}

\bd{Keywords} \quad polynomial roots, bounds on zeros, power series, polynomial reciprocals, analytic combinatorics, meromorphic functions\par

\section{Introduction}
We briefly review the bounds proved by \citet{kalantari_2004, kalantari_2008} on the moduli of polynomial roots. \par 
Let $p(z) = a_nz^n + \cdots + a_1z + a_0,\; a_na_0\neq 0$, be an arbitrary polynomial of degree $n$ and $\rho$ be an arbitrary zero of $p(z)$.  Moreover, let $r_m$ be, as in \citet{hohertz_kalantari}, the unique positive root of $f_m(x):=x^{m+1} + x - 1$, a root that increases to 1 as $m\rightarrow\infty$ and necessarily lies on the interval $[\frac{1}{2}, 1)$.  Finally, let
	\begin{equation} 
		D_{m}(z) := \det
		\begin{pmatrix}
			p'(z) & \frac{p''(z)}{2!} & \cdots & \frac{p^{(m-1)}(z)}{(m-1)!} & \frac{p^{(m)}(z)}{m!} \\
			p(z) & p'(z) & \cdots & \ddots & \frac{p^{(m-1)}(z)}{(m-1)!} \\
			0 & p(z) & \ddots & \ddots & \vdots \\
			\vdots & \vdots & \ddots & \ddots & \frac{p''(z)}{2!} \\
			0 & 0 & \cdots & p(z) & p'(z)
		\end{pmatrix}
	\end{equation}
	and
	\begin{equation} 
		\widehat{D}_{m,j}(z) := \det
		\begin{pmatrix}
			\frac{p''(z)}{2!} & \frac{p'''(z)}{3!} & \cdots & \frac{p^{(m)}(z)}{m!} & \frac{p^{(j)}(z)}{j!} \\
			p'(z) & \frac{p''(z)}{2!} & \ddots & \frac{p^{(m-1)}(z)}{(m-1)!} & \frac{p^{(j-1)}(z)}{(j-1)!} \\
			p(z) & p'(z) & \ddots & \vdots & \vdots \\
			\vdots & \vdots & \ddots & \frac{p''(z)}{2!} & \frac{p^{(j-m+2)}(z)}{(j-m+2)!}\\
			0 & 0 & \cdots & p'(z) & \frac{p^{(j-m+1)}(z)}{(j-m+1)!}
		\end{pmatrix}.
	\end{equation}
Kalantari (see \citet{kalantari_2004}, \citet{kalantari_2008}, and \citet{jin}) proves that
	\begin{equation} \label{eq:Lm}
		\aval{\rho} \geq \frac{r_m}{\gamma_m} =: L_m
	\end{equation}
for $m\geq 2$, where
	\begin{equation}
		\gamma_m:=\max_{k=m+1, \cdots, m+n} \aval{\frac{\widehat{D}_{m+1,k+1}\inps{0}}{D_{m+1}\inps{0}}}^{\frac{1}{k}}.\footnote{Clearly, one obtains a dual upper bound $U_m$ on $\aval{\rho}$ by applying the previous formulas to the reciprocal polynomial of $p(z)$ and reciprocating both sides of the inequality of Equation \eqref{eq:Lm}.}
	\end{equation}
Using the equation
	\begin{equation} \label{eq:convo}
		\inps{-1}^{m}\frac{\widehat{D}_{m+1,k+1}\inps{0}}{D_{m+1}\inps{0}}= \sum_{j=\max\sbld{0,k-n}}^m a_{k-j}b_j,
	\end{equation}
where $b_j$ is the coefficient of $z^j$ in the Maclaurin series of $\frac{1}{p(z)}$, \citet{jin} proves that $L_m$ converges, as $m\rightarrow\infty$, to $\rhomin:=\min_{\rho : p(\rho) = 0} \aval{\rho}$ ,\footnote{\ic{resp.}, $U_m$ to $\rhomax:=\max_{\rho : p(\rho) = 0} \aval{\rho}$}  and thus is asymptotically sharp.  In effect, he proves the following theorem:
	\begin{theorem}[equivalent to Theorem 3.1, \citet{jin}] \label{thm:jin}
		For every polynomial $p(z)$ and positive $\epsilon$ close to zero, there exists a positive number $m_{\epsilon}$ such that $m\geq m_{\epsilon}$ implies
		\begin{equation}
			\frac{L_m}{\rhomin} \geq 1 - \epsilon.
		\end{equation}
	\end{theorem}
Jin's proof, however, does not provide an algorithm for finding $m_{\epsilon}$.  In this paper we provide such an algorithm in the form of an equation, thus proving the following theorem:
	\begin{theorem} \label{thm:main}
		For every positive integer $n$ and triple $\inps{\epsilon, \alpha, \beta}$ of positive reals with $\epsilon$ close to zero, there exists a positive number $m_{\epsilon}$ such that $m\geq m_{\epsilon}$ implies
		\begin{equation} \label{eq:thmmain}
			\frac{L_m}{\rhomin} \geq 1 - \epsilon
		\end{equation}
	for any polynomial $p(z)$ with degree $n$, $\alpha:=\max\sbld{\aval{a_0}, \cdots, \aval{a_n}}$, and $\beta:=\max\sbld{1,\rhomin^n}\cdot \sum_{j=0}^{n-1} \rhomin^j$.
	\end{theorem}
Note that this theorem implies Theorem \ref{thm:jin}.  \par
We acknowledge a few drawbacks of our formula: in particular, the value of $m_{\epsilon}$ it provides, though sufficient, does not appear to be optimal.  We therefore devote a section of this paper to experimental results on optimal $m_{\epsilon}$; in particular, we conjecture that $m_{\epsilon} = O\inps{\frac{1}{\epsilon^d}}$ for some $d\ll 2$.  Since Jin's method for calculating $L_m$ runs in $O(mn)$ time, the truth of our conjecture would imply that roots of degree $n$ polynomials could be bounded with error tolerance $\epsilon$ in $O\inps{\frac{n}{\epsilon^d}}$ time. \par

(In this paper, $p(z)$, $n$, $\alpha$, \ic{etc.}, retain the definitions they are assigned in this introduction unless otherwise specified.) 

\section{Main results}
From here on, we assume $a_0 = p(0) = 1$ and $n > 2$ without loss of generality.\par
By the proof of Theorem 3.1 of \citet{jin}, 
	\begin{equation}
		\gamma_m \leq \max_{k=1+m,\cdots,n+m} \left[\alpha\beta Q(k) \right]^{1/k}\cdot \rhomin^{-1},
	\end{equation}
where 
	\begin{align}
		\alpha &:= \max\sbld{1, \aval{a_1}, \cdots, \aval{a_n}}, \\
		\beta&:=\max\sbld{1,\rhomin^n}\cdot \sum_{j=0}^{n-1} \rhomin^j,
	\end{align}
and $Q(k)$ is a monotonically increasing function\footnote{in \citet{jin}, $Q(k)$ is a polynomial with positive coefficients, of degree one less than the maximum multiplicity of the roots of $p(z)$.  However, the proof of his Theorem 3.1 requires only that $Q(k)$ be increasing.} such that
	\begin{equation}
		\aval{b_k}\leq Q(k)\cdot \rhomin^{-k}.
	\end{equation}
	
	\begin{lemma} \label{lemma:bndms}
		Let $d_k(n):=\binom{n+k-1}{k}$ be the number of $k$-multisets of members of $\sbld{1, \dots, n}$.  Then
		\begin{equation}
			d_k(n) < \frac{\inps{k+\frac{n}{2}}^{n-1}}{\Gamma(n)}.
		\end{equation}
		\begin{proof}
			This is the first case of Lemma 2 of \citet{grinshpan}.
		\end{proof}
	\end{lemma}
	
	\begin{theorem} \label{thm:bndrec}
		Let $p(z)$ be the degree $n$ polynomial \label{thm:recip}
		\begin{equation}
			p(z):=a_n\prod_{i=1}^n \inps{z - z_i} = a_nz^n + a_{n-1}z^{n-1} + \cdots + 1,
		\end{equation}
		and let $\rhomin$ be the least of the moduli of the roots of $p(z)$.  Moreover, let
		\begin{equation}
			\frac{1}{p(z)} = 1 + b_1z + \cdots
		\end{equation}
	be the Maclaurin series of $\frac{1}{p(z)}$.  Then
		\begin{equation}
			|b_k|\leq \rhomin^{-k}\cdot\frac{\inps{k+\frac{n}{2}}^{n-1}}{\Gamma\inps{n}}.
		\end{equation}
		\begin{proof}
 		Clearly,
			\begin{align}
				b_k &= \frac{1}{k!}\cdot\left[\frac{d^k}{dz^k}\inps{\frac{1}{p(z)}}\right]_{z = 0}.
			\end{align}
		By the general Leibniz Product Rule,
		\begin{align}	
			\frac{d^k}{dz^k}\inps{\frac{1}{p(z)}} &= \frac{1}{a_n}\cdot \sum_{k_1 + \cdots + k_n = k} \binom{k}{k_1, \cdots, k_n} \prod_{i=1}^n\frac{d^{k_i}}{dz^{k_i}}\inps{z - z_i}^{-1} \\
			&=\frac{1}{a_n} \cdot\sum_{k_1 + \cdots + k_n = k} \binom{k}{k_1, \cdots, k_n} \prod_{i=1}^n\frac{(-1)^{k_i}k_i !}{(z - z_i)^{k_i + 1}}  \\
			&=\frac{ \inps{-1}^k k!}{a_n}\cdot  \sum_{k_1 + \cdots + k_n = k} \prod_{i=1}^n\frac{1}{(z - z_i)^{k_i + 1}} \\
			&=\frac{ \inps{-1}^k k!}{p(z)}\cdot  \sum_{k_1 + \cdots + k_n = k} \prod_{i=1}^n\frac{1}{(z - z_i)^{k_i}}.
			\intertext{Dividing both sides by $k!$, setting $z = 0$, and taking absolute values and applying the Triangle Inequality, we obtain}
			|b_k| &\leq \sum_{k_1 + \cdots + k_n = k} \prod_{i=1}^n \frac{1}{|z_i|^{k_i}} \\
			&\leq  \sum_{k_1 + \cdots + k_n = k} \rhomin^{-k}.
			\intertext{Since the sum is taken over all size $k$ multisets of $\sbld{1,\cdots,n}$, this last inequality reduces to}
			|b_k|  &\leq  \binom{n+k-1}{k} \rhomin^{-k},
			\intertext{from which it follows, by Lemma \ref{lemma:bndms}, that}
			|b_k|  &\leq \frac{\inps{k+\frac{n}{2}}^{n-1}}{\Gamma(n)}\cdot \rhomin^{-k}.
		\end{align}
		\end{proof}
	\end{theorem}
By Theorem \ref{thm:bndrec}, we may set 
	\begin{equation}
		Q(k):=\frac{\inps{k+\frac{n}{2}}^{n-1}}{\Gamma\inps{n}},
	\end{equation} 
which, as required, is manifestly an increasing function of $k$.
We therefore seek integral $k$ from $1+m$ to $n+m$ that maximizes 
	\begin{equation}\label{eq:rt}
		\left[\alpha\beta\cdot\frac{\inps{k+\frac{n}{2}}^{n-1}}{\Gamma(n)} \right]^{1/k}\cdot \rhomin^{-1}. 
	\end{equation}
\begin{lemma} \label{lemma:dec}
	Let $c_1, c_2, c_3 > 0$.  Then, for large $m$, the function
	\begin{equation}
		f(m):=\left[c_1\cdot \inps{m + c_2}^{c_3}\right]^{1/m}
	\end{equation}
	is decreasing.
	\begin{proof}
	We take the logarithm of $f(m)$, multiply both sides by $m$, and derive:
	\begin{align}
		\log\inps{f(m)} &= \frac{\log\inps{c_1} + {c_3}\cdot\log\inps{m+c_2}}{m} \\
		m\log\inps{f(m)} &= \log\inps{c_1} + {c_3}\cdot\log\inps{m+c_2}
	\end{align}
	\begin{equation}
		\frac{m\cdot f'(m)}{f(m)} + \log\inps{f(m)}= \frac{{c_3}}{m+c_2}.
	\end{equation}
	Thus
	\begin{align}
		\frac{m\cdot f'(m)}{f(m)} &= \frac{{c_3}}{m+c_2} -  \frac{\log\inps{c_1} + {c_3}\cdot\log\inps{m+c_2}}{m} \\
			f'(m) &= \frac{f(m)}{m}\cdot\inps{\frac{{c_3}}{m+c_2} -  \frac{\log\inps{c_1} + {c_3}\cdot\log\inps{m+c_2}}{m}}.
	\end{align}
	Asymptotically, the term $-c_3\cdot\log\inps{m+c_2}$ dominates, so that $f'(m)$ is negative for large $m$.
	\end{proof}
\end{lemma}

\begin{lemma} \label{lemma:rm}
	For large $m$,
	\begin{equation}
		r_m \geq \inps{m+2}^{-\frac{1}{m+1}}.
	\end{equation}
	\begin{proof}
		Since $f_m\inps{x}$ is increasing for small positive $x$, it suffices to show that 
		\begin{equation}
			\inps{\inps{m+2}^{-\frac{1}{m+1}}}^{m+1} + \inps{m+2}^{-\frac{1}{m+1}} - 1 \leq 0.
		\end{equation}
		Setting $u(t):= \inps{t+2}^{-\frac{1}{t+1}}$ and $f\inps{u}:=u^{t+1} + u - 1$, it suffices to prove that
		\begin{enumerate}
			\item $\frac{df}{dt} > 0$ and
			\item $\lim_{t\rightarrow\infty} f(t) = 0$.
		\end{enumerate}
		Now,
		\begin{equation}
			\frac{df}{dt} = \frac{\partial f}{\partial u}\cdot \frac{\partial u}{\partial t} + \frac{\partial f}{\partial t},
		\end{equation}
		where
		\begin{align}
			\frac{\partial f}{\partial u} &= \inps{t+1}u^t + 1 \\
			 &= \inps{t+1}\inps{t+2}^{-\frac{t}{t+1}} + 1
		\intertext{and}
			\frac{\partial f}{\partial t} &= u^{t+1}\cdot\logf{u} \\
			&= \frac{1}{t+2} \cdot -\frac{\logf{t+2}}{t+1}.
		\end{align}
		To calculate $\frac{\partial u}{\partial t}$, we take the logarithm of both sides of the equation
		\begin{equation}
			u = \inps{t+2}^{-\frac{1}{t+1}},
		\end{equation}
		obtaining
		\begin{align}
			\log u &= -\frac{\logf{t+2}}{t+1} \\
			\frac{1}{u}\cdot \frac{\partial u}{\partial t} &= - \frac{\inps{t+1}\cdot \frac{1}{t+2} - \logf{t+2}}{\inps{t+1}^2}.
		\end{align}
		As $t\rightarrow\infty$, $\frac{\partial f}{\partial u}\rightarrow 2$, $\frac{\partial u}{\partial t}$ is positive, and $\frac{\partial f}{\partial t} \rightarrow 0$.  Thus $\frac{df}{dt} > 0$ for sufficiently large $t$. \par
		As for the second prong,
		\begin{align}
			f\inps{t} &= \frac{1}{t+2} + \inps{t+2}^{-\frac{1}{t+1}} - 1 \\
			&\rightarrow 0 + 1 - 1 = 0.
		\end{align}
	\end{proof}
\end{lemma}

Suppose $m$ sufficiently large that the conclusions of Lemmas \ref{lemma:dec} and \ref{lemma:rm} apply.  By Lemma \ref{lemma:dec}, the quantity of Equation \eqref{eq:rt} is maximized for $k = 1+m$.  Thus

\begin{equation}	
	\gamma_m \leq \left[\alpha\beta\frac{\inps{m+1+\frac{n}{2}}^{n-1}}{\Gamma\inps{n}}\right]^{\frac{1}{m+1}}\cdot \rhomin^{-1},
\end{equation} implying, by Lemma \ref{lemma:rm},  that

\begin{align}	
	L_m:=\frac{r_m}{\gamma_m} &\geq \left[\frac{\Gamma\inps{n}}{\alpha\beta\inps{m+2}\inps{m+1+\frac{n}{2}}^{n-1}}\right]^{\frac{1}{m+1}}\cdot \rhomin \\
	\frac{L_m}{\rhomin}&\geq \left[\frac{\Gamma\inps{n}}{\alpha\beta\inps{m+1+\frac{n}{2}}^{n}}\right]^{\frac{1}{m+1}}=: g(m;n,\alpha,\beta).\label{eq:eps} 
\end{align}

Therefore, Equation \eqref{eq:thmmain} of Theorem \ref{thm:main} holds provided that $m\geq \max\sbld{m_{\epsilon},m_{\ell}}$, where $m_{\epsilon}$ solves the equation
\begin{equation} \label{eq:meps}
	g(m_{\epsilon};n,\alpha,\beta) = 1-\epsilon
\end{equation} 
and $m_{\ell}$ is the least $m$ that is ``sufficiently large" for the conclusions of Lemmas \ref{lemma:dec} and \ref{lemma:rm} to hold.  Since $g(m;n,\alpha,\beta)$ is the reciprocal of the $m^{th}$ root of a polynomial in $m$, it has a limit of 1 as $m\rightarrow\infty$; and since both $m_{\epsilon}$ and $m_{\ell}$ depend only on $n$, $\alpha$, $\beta$, and $\epsilon$, Theorem \ref{thm:main} is proved.  (We hereafter assume $m_{\epsilon}\geq m_{\ell}$.)  \par

\section{Discussion}
 We anticipate three primary objections to our method:
\begin{enumerate}
	\item \ic{Equation \eqref{eq:meps} does not yield an obvious asymptotic bound for $m_{\epsilon}$.}  Indeed, we have not yet succeeded in rigorously proving such a bound.  However, the experimental results we detail in Section \ref{sec:expm} suggest that $m_{\epsilon} = O\inps{\frac{C}{\epsilon^d}}$ for some constant $C\leq 1.1$ and $d\in\inps{1.28,1.42}$.  Note that this estimate does not depend on $n$.
	\item \ic{The function $g(m_{\epsilon};n,\alpha,\beta)$ approximates $\frac{L_{m_{\epsilon}}}{\rhomin}$ yet depends on $\beta$, which is itself a function of $\rhomin$.} This apparent circularity can be overcome by replacing $\rhomin$ in the formula for $\beta$ with the modulus of a known root, or a known upper bound of $\rhomin$.  In particular, this allows the estimate $\beta = n$ if at least one root of $p(z)$ is known to lie in $\overline{\D}$. \par
	Consider, for example, $p(z) = 3z^7 - z^2 + 2$.  By Viete's formulas, $p(z)$ has at least one zero in $\D$; thus we may estimate $\beta$ with $n = 7$ to obtain the equation
		\begin{equation}
			\left[\frac{\Gamma\inps{7}}{\frac{3}{2}\cdot 7\cdot \inps{m_{\epsilon}+1+\frac{7}{2}}^{7}}\right]^{\frac{1}{m_{\epsilon}+1}} = 1 - \epsilon
		\end{equation}		
	for $m_{\epsilon}$.  For $\epsilon = 0.05$, the solution is $m_{\epsilon} = 828$; indeed, $L_{828} \approx 0.876282$ and $\rhomin \approx 0.88$, yielding 
		\begin{equation}
			\frac{L_{828}}{\rhomin} \approx 0.995817 \gg 0.95.
		\end{equation}
	\item \ic{The function $g(m_{\epsilon};n,\alpha,\beta)$ tends to overestimate the minimum $m$ necessary for error tolerance $\epsilon$.}  Returning to the previous example, for $p(z) = 3z^7 - z^2 + 2$ only $m=20$ calculations, and not $m_{\epsilon} = 828$, are required to approximate $\rhomin \approx 0.88$ to $95\%$ accuracy.  In the next section we suggest a formula for an $m$ (generally less than $m_{\epsilon}$) that tends to suffice in practice.
\end{enumerate}

\section{Experimental results} \label{sec:expm}
\begin{conj}\label{conj}
	Equation \eqref{eq:thmmain} holds for $m = m_c$, where
	\begin{equation}
		m_{c} = O\inps{\frac{1}{\epsilon^c}}
	\end{equation}
	for some fixed $1 < c \ll 2$.  In particular, $m_c$ does not depend on $n$.
\end{conj}
To arrive at Conjecture \ref{conj}, we performed two experiments:
\begin{enumerate}
	\item Took the average value of $m_c$ for 100 polynomials with random degree between 2 and 256, constant term 1, and non-constant coefficients of random maximum absolute value between 2 and 256, setting for each polynomial $\epsilon = 0.2\cdot 2^{-j}$ for some random integer $j$ between 0 and 7 inclusive. \label{item:all}
	\item For each $j = 0, ..., 19$, took the average value of $m_c$ for 200 polynomials with degree 10 and coefficients of absolute value no more than $99$, setting $\epsilon = 0.2\cdot \inps{0.8}^j$. \label{item:eps}
\end{enumerate}

In Experiment \ref{item:all} we performed linear regression on $m_c$ with independent variables $\alpha$, $n$, and $\epsilon$; our models and the corresponding values of $r^2$ are shown in Table \ref{tab:all}.  Most strikingly, the results of Experiment \ref{item:all} suggest no correlation between $m_c$ and the degree of polynomial, attributing less than 1\% of variation in $\log m_c$ to changes in variables other than $\epsilon$.  Using these results, we might infer the following approximate formula for $m_c$ (ignoring variables for degree and $\alpha$, whose correlation coefficients are small):
\begin{equation} \label{eq:expm}
	m_c = \frac{C}{\epsilon^d},
\end{equation}
where $C\in (1.071929, 1.095269)$ and $d\in(1.286622, 1.287167)$.

By contrast, linear regression on the results of Experiment \ref{item:eps} yields the formula of Equation \eqref{eq:expm} with $C\approx 0.540619$ and $d \approx 1.412172$ ($r^2 = 0.9961$).  Taking $C$ from Experiment \ref{item:all} and $d$ from Experiment \ref{item:eps}, we define
\begin{equation}
	m_{exp} := \frac{1.095269}{\epsilon^{1.412172}}
\end{equation} as an approximation of suitable $m$ for purposes of the next section.

\section{Examples} 
In \citet{hohertz_kalantari} we introduced the {Collatz polynomials} $P_N(z)$, defined\footnote{that is, defined at least for those positive integers for which the Collatz trajectory eventually terminates - conjectured to be all of them.} as
\begin{equation}
	P_N(z) := \sum_{j=0}^{h(N)} c^j\inps{N}\cdot z^j
\end{equation}
where 
\begin{equation}
c(N) :=
\begin{cases}
	\frac{3N+1}{2}, & N\mbox{ odd} \\
	\frac{N}{2}, & N\mbox{ even,}
\end{cases}
\end{equation} 
$c^{j+1}\inps{N}:=c\inps{c^j\inps{N}}$, and $h(N):=\min \sbld{j\::\:c^{j}\inps{N}=1}$.  Defining $M(N):=\max_{j\geq 0} c^j\inps{N}$, we have that
\begin{align}
	\alpha &= \frac{M(N)}{N}, \\
	\beta &\approx 2^{h(N)},
\intertext{and}
	n &= h(N).
\end{align} Using the expected value $h(N)=\frac{2}{\logf{\frac{4}{3}}}\cdot \log N \approx 6.952118\cdot \log N$ proposed in \citet{lagarias} and the value $\alpha = 8N$ conjectured in \citet{silva} to bound $\frac{M(N)}{N}$ above, we obtain the formula
\begin{equation} \label{eq:collatzeq}
	 \left[\frac{\Gamma\inps{6.952118\cdot\log N}}{8N\cdot \inps{2m_{\epsilon}+2+6.952118\cdot\log N}^{6.952118\cdot\log N }}\right]^{\frac{1}{m_{\epsilon}+1}} = 1 - \epsilon.
\end{equation}
For $P_5\inps{z} = 5 + 8z + 4z^2 + 2z^3 + z^4$ and $\epsilon = 0.05$, Equation \eqref{eq:collatzeq} has solution $m_{\epsilon} = 1518$ (rounded up to the nearest integer).  This value gives the estimate $L_{1434} \approx 0.995717$, within $0.5\%$ of $\rhomin = \aval{-1} = 1$.  On the other hand, $\lceil m_{exp} \rceil = 76$; the corresponding estimate is $L_{76} \approx 0.947933$, just over $5\%$ less than $\rhomin$. \par

Consistent with Conjecture \ref{conj} and the results of Section \ref{sec:expm}, $m_{exp}$ appears to bound $\rhomin$ well for the polynomials $P_N$.  Indeed, on the interval $N\in\left[2,703\right]$, the degree of $P_N$ ranges from 1 to 108, yet $\frac{L_{m_{exp}}}{\rhomin} < 0.95$ for only six values of $N$ (it performs worst for $N = 137$, for which $\frac{L_{m_{exp}}}{\rhomin} \approx 0.927254$).

\section{Conclusion}
All of our experiments suggest the independence of the value $m$ in Jin's method from the polynomial degree $n$: this striking conjecture would have powerful implications if true and warrants further study.  In particular, we encourage research into proof or disproof of Conjecture \ref{conj}, as well as an asymptotic bound on the minimum $m$ required for error tolerance $\epsilon$ (or a value for the exponent $d$ if Conjecture \ref{conj} holds).  Finally, acknowledging the theoretical significance of our Equation \eqref{eq:meps}, we nevertheless hope to sharpen the resulting $m_{\epsilon}$, for which we would like to find a closed formula.

\section{Tables}
\begin{table}[htp]
\caption{Results of Experiment \ref{item:all}.  Average $m_c$ for polynomials with random values of $\alpha$, $n$, and $\epsilon$.}
\begin{center}
\begin{tabular}{|c|c|c|c|c|}
	\hline
	\multicolumn{3}{|c|}{\textsc{variable meanings}} & \multicolumn{2}{|c|}{\textsc{linear models}} \\ \hline 
	$u$ & $v$ & $w$ & $\log m_c =$ & $r^2$ \\ \hline
	$\alpha$ & - & - & $6.0405 - 0.0025u$ & 0.0063 \\ \hline
	- & \textsc{degree} & - & $5.9261 - 0.0015v$ & 0.0028 \\ \hline
	- & - & $\log\epsilon$ & $0.0910 - 1.2866w$ & 0.9958 \\ \hline
	$\alpha$ & \textsc{degree} & - & $6.2239 - 0.0025u - 0.0015v$ & 0.0090 \\ \hline
	$\alpha$ & - & $\log\epsilon$ & $0.0728 + 0.0001 u - 1.2871w$ & 0.9958 \\ \hline
	 - & \textsc{degree} & $\log\epsilon$ & $0.0886 + 0.00002v - 1.2867w$ & 0.9958 \\ \hline
	  $\alpha$ & \textsc{degree} & $\log\epsilon$ & $0.0695 + 0.0001u + 0.00002v - 1.2872w$ & 0.9958 \\ \hline
\end{tabular}
\end{center}
\label{tab:all}
\end{table}

\begin{table}[htp]
\caption{Results of Experiment \ref{item:eps}.  Average $m_c$ for degree 10 polynomials with coefficients $\aval{a_i} \leq 99$ and $\epsilon = 0.2\cdot\inps{0.8}^j$.}
\begin{center}
\begin{tabular}{|c|c||c|c||c|c||c|c|}
\hline 
	j & ${m_c}_{avg}$ & j & ${m_c}_{avg}$ & j & ${m_c}_{avg}$ & j & ${m_c}_{avg}$ \\ \hline
	0 & 4.085 & 5 & 27.525 & 10 & 138.140  & 15 & 586.480  \\ \hline 
	1 & 6.090 & 6 & 38.270 & 11 & 186.495 & 16 & 773.515  \\ \hline 
	2 & 9.105 & 7 & 53.895 & 12 & 250.010 & 17 & 1018.430  \\ \hline 
	3 & 13.125  & 8 & 74.775 & 13 & 333.670  & 18 & 1336.235  \\ \hline
	4 & 19.050 & 9 & 101.935 & 14 & 443.225 & 19 & 1748.345  \\ \hline 
\end{tabular}
\end{center}
\label{tab:eps}
\end{table}

\clearpage
\setcitestyle{numbers}
\bibliography{JinMethod-Oct06}

\begin{thebibliography}{7}
\providecommand{\natexlab}[1]{#1}
\providecommand{\url}[1]{\texttt{#1}}
\expandafter\ifx\csname urlstyle\endcsname\relax
  \providecommand{\doi}[1]{doi: #1}\else
  \providecommand{\doi}{doi: \begingroup \urlstyle{rm}\Url}\fi

\bibitem[Kalantari(2004)]{kalantari_2004}
Bahman Kalantari.
\newblock An infinite family of bounds on zeros of analytic functions and
  relationship to {S}male's bound.
\newblock \emph{Mathematics of Computation}, 74\penalty0 (250):\penalty0
  841–853, 2004.
\newblock \doi{10.1090/s0025-5718-04-01686-2}.

\bibitem[Kalantari(2008)]{kalantari_2008}
Bahman Kalantari.
\newblock \emph{Polynomial {R}oot-finding and {P}olynomiography}.
\newblock World Scientific, Hackensack, NJ, 2008.

\bibitem[Hohertz and Kalantari(2020)]{hohertz_kalantari}
Matt Hohertz and Bahman Kalantari.
\newblock Collatz polynomials: an introduction with bounds on their zeros.
\newblock 2020.

\bibitem[Jin(2006)]{jin}
Yi~Jin.
\newblock On efficient computation and asymptotic sharpness of {K}alantari's
  bounds for zeros of polynomials.
\newblock \emph{Mathematics of Computation}, 75:\penalty0 1905–1912, 2006.

\bibitem[Grinshpan(2010)]{grinshpan}
Arcadii~Z. Grinshpan.
\newblock Weighted inequalities and negative binomials.
\newblock \emph{Advances in Applied Mathematics}, 45:\penalty0 564–606, 2010.

\bibitem[Applegate and Lagarias(2002)]{lagarias}
David Applegate and Jeffrey~C. Lagarias.
\newblock Lower bounds for the total stopping time of $3x + 1$ iterates.
\newblock \emph{Mathematics of Computation}, 72\penalty0 (242):\penalty0
  1035–1050, 2002.
\newblock \doi{10.1090/s0025-5718-02-01425-4}.

\bibitem[Silva(1999)]{silva}
Tomás Oliveira~E Silva.
\newblock Maximum excursion and stopping time record-holders for the $3x+1$
  problem: Computational results.
\newblock \emph{Mathematics of Computation}, 68\penalty0 (225):\penalty0
  371–385, 1999.
\newblock \doi{10.1090/s0025-5718-99-01031-5}.

\end{thebibliography}

\end{document}